\documentclass{amsart}

\def\uu{\mathcal U}
\def\uuu{\overline{\mathcal U}}

\newcommand{\remove}[1]{ }

\newtheorem{theorem}{Theorem}[section]

\newtheorem{lemma}[theorem]{Lemma}

\theoremstyle{definition}

\theoremstyle{remark}
\newtheorem*{remark}{Remark}
\newtheorem*{remarks}{Remarks}

\numberwithin{equation}{section}

\begin{document}
\title{A property of algebraic univoque numbers}
\author{Martijn de Vries}
\address{Delft University of Technology, Mekelweg 4, 2628 CD Delft, the Netherlands}
\email{w.m.devries@ewi.tudelft.nl}
\date{\today}
\thanks{}

\begin{abstract}
Consider the set $\uu$ of real numbers $q \ge 1$ for which only one sequence  $(c_i)$ of integers $0 \le c_i \le q$ satisfies the equality $\sum_{i=1}^{\infty} c_i q^{-i} = 1$. In this note we show that the set of algebraic numbers in $\uu$ is dense in the closure $\uuu$ of $\uu$.
\end{abstract}

\maketitle


\section{Introduction}\label{s1}

Given a real number $q \ge 1$, a $q-${\it expansion} (or simply {\it expansion}) is a sequence $(c_i)=c_1 c_2 \ldots$ of integers satisfying $0 \le c_i \le q$ for all $i \geq 1$ such that 
\begin{equation*}
\frac{c_1}{q} + \frac{c_2}{q^2} + \frac{c_3}{q^3} + \cdots = 1.
\end{equation*}

One such expansion, denoted by $(\gamma_i(q))= (\gamma_i)$, is obtained by performing the {\it greedy algorithm} of R\'enyi 
(\cite{R}): if $\gamma_i$ is already defined for $ i < n$, then $\gamma_n$ is the largest integer satisfying
\begin{equation*}
\sum_{i=1}^{n} \frac{\gamma_i}{q^i} \leq 1.
\end{equation*}
Equivalently, $(\gamma_i)$ is the largest expansion in lexicographical order. 

If $q >1$, then another such expansion, denoted by $(\alpha_i(q))= (\alpha_i)$, is obtained by performing the {\it quasi-greedy algorithm}: if $\alpha_i$ is already defined for $i < n$, then $\alpha_n$ is the largest integer satisfying 
\begin{equation*}
\sum_{i=1}^{n} \frac{\alpha_i}{q^i} < 1.
\end{equation*}

An expansion is called {\it infinite} if it contains infinitely many nonzero terms; otherwise it is called {\it finite}. Observe that there are no infinite expansions if $q = 1$: the only 1-expansions are given by 
$10^{\infty}, 010^{\infty}, 0010^{\infty}, \ldots$. 
On the other hand, if $q >1$, then $(\alpha_i)$ is the largest infinite expansion in lexicographical order. 

For any given $q >1$, the following relations between the quasi-greedy expansion and the greedy expansion are straightforward. The greedy expansion is finite if and only if $(\alpha_i)$ is periodic. If $(\gamma_i)$ is finite and $\gamma_m$ is its last nonzero term, then $m$ is the smallest period of $(\alpha_i)$, and
\begin{equation*}
\alpha_i = \gamma_i \quad \mbox{for } i = 1, \ldots, m-1, \quad \mbox{and } \alpha_m = \gamma_m -1.
\end{equation*}

Erd\H{o}s, Horv\'ath and Jo\'o (\cite{EHJ}) discovered that for some real numbers $q >1 $ there exists only one $q-$expansion. Subsequently, the set $\uu$ of such {\it univoque numbers} was characterized in \cite{EJK1}, \cite{EJK2}, \cite{KL3} (see Theorem~\ref{t21}). Using this characterization, Komornik and Loreti showed in \cite{KL1} that $\uu$ has a smallest element $q' \approx 1.787$ and the corresponding expansion $(\tau_i)$ 
is given by the truncated Thue-Morse sequence, defined by setting $\tau_{2^N}=1$ for $N=0,1,\ldots$ and 
\begin{equation*}
\tau_{2^N + i} = 1 - \tau_i \quad \mbox{for }1 \le i < 2^N, \, N=1,2, \ldots.
\end{equation*}
Allouche and Cosnard (\cite{AC}) proved that the number $q'$ is transcendental. This raised the question whether there exists a smallest algebraic univoque number. Komornik, Loreti and Peth\H{o} (\cite{KL2}) answered this question in the negative by constructing a decreasing sequence $(q_n)$ of algebraic univoque numbers converging to $q'$. 

It is the aim of this note to show that for each $q \in \uu$ there exists a sequence of algebraic univoque numbers converging to $q$:

\begin{theorem}\label{t11}
The set $\mathcal{A}$ consisting of all algebraic univoque numbers is dense in $\uuu$.
\end{theorem}

Our proof of Theorem \ref{t11} relies on a characterization of the closure $\uuu$ of $\uu$, recently obtained by Komornik and Loreti in \cite{KL3} (see Theorem~\ref{t22}). 

\section{Proof of Theorem~\ref{t11}}\label{s2}

In the sequel, a sequence always means a sequence of nonnegative integers. We use systematically the lexicographical order between sequences; we write $(a_i) < (b_i)$ if there exists an index $n \geq 1$ such that $a_i=b_i$ for $i < n$ and $a_n < b_n$. This definition extends in the obvious way to sequences of finite length. 

The following algebraic characterization of the set $\uu$ can be found in \cite{EJK1}, \cite{EJK2}, \cite{KL3}:

\begin{theorem}\label{t21}
The map $q \mapsto (\gamma_i(q))$ is a strictly increasing bijection between the set $\uu$ and the set of all sequences $(\gamma_i)$ satisfying 
\begin{equation}\label{21}
\gamma_{j+1} \gamma_{j+2} \ldots < \gamma_1 \gamma_2 \ldots \quad \mbox{for all }\, j \geq 1
\end{equation}
and 
\begin{equation}\label{22}
\overline{\gamma_{j+1} \gamma_{j+2} \ldots} < \gamma_1 \gamma_2 \ldots \quad \mbox{for all }\, j \geq 1
\end{equation}
where we use the notation $\overline{\gamma_n} := \gamma_1 - \gamma_n$. 
\end{theorem}

\begin{remark} 
It was essentially shown by Parry (see \cite{P}) that a sequence $(\gamma_i)$ is the greedy $q$-expansion for some $q \ge 1$ if and only if $(\gamma_i)$ satisfies the condition \eqref{21}. 
\end{remark}

Using the above result, Komornik and Loreti (\cite{KL3}) investigated the topological structure of the set $\uu$. In particular they showed that $\uuu \setminus \uu$ is dense in $\uuu$. Hence the set $\uuu$ is a perfect set. Moreover, they established an analogous characterization of the closure $\uuu$ of $\uu$:
\begin{theorem}\label{t22}
The map $q \mapsto (\alpha_i(q))$ is a strictly increasing bijection between the set $\uuu$ and the set of all sequences $(\alpha_i)$ satisfying 
\begin{equation}\label{23}
\alpha_{j+1} \alpha_{j+2} \ldots \le \alpha_1 \alpha_2 \ldots \quad \mbox{for all }\, j \geq 1
\end{equation}
and 
\begin{equation}\label{24}
\overline{\alpha_{j+1} \alpha_{j+2} \ldots} < \alpha_1 \alpha_2 \ldots \quad \mbox{for all }\, j \geq 1
\end{equation}
where we use the notation $\overline{\alpha_n} := \alpha_1 - \alpha_n$. 
\end{theorem}

\begin{remarks}
\mbox{}
\begin{itemize}
 
\item It was shown in \cite{BK} that a sequence $(\alpha_i)$ is the quasi-greedy $q$-expansion for some $q >1$ if and only if $(\alpha_i)$ is infinite and satisfies \eqref{23}. Note also that a sequence satisfying \eqref{23} and \eqref{24} is automatically infinite.  

\item If $q \in \uuu \setminus \uu$, then we must have equality in \eqref{23} for some $j \geq 1$, i.e., the greedy $q$-expansion is finite for each $q \in \uuu \setminus \uu$. On the other hand, it follows from Theorems~\ref{t21} and \ref{t22} that a sequence of the form $(1^n0)^{\infty}$ $(n \ge 2)$ is the quasi-greedy $q$-expansion for some $q \in \uuu \setminus \uu$. Hence the set $\uuu \setminus \uu$ is countably infinite.  
\end{itemize}

\end{remarks} 

The following technical lemma is a direct consequence of Theorem~\ref{t22} and Lemmas 3.4 and 4.1 in \cite{KL3}:

\begin{lemma}\label{l23}
Let $(\alpha_i)$ be a sequence satisfying \eqref{23} and \eqref{24}. Then 
\begin{itemize}
\item[\rm (i)] there exist arbitrary large integers $m$ such that 
\begin{equation}\label{25}
\overline{\alpha_{j+1} \ldots \alpha_m} < \alpha_1 \ldots \alpha_{m-j} \quad \mbox{for all }\, 0 \le j < m; 
\end{equation}
\item[\rm (ii)] for all positive integers $m \geq 1$, 
\begin{equation}\label{26}
\overline{\alpha_{1} \ldots \alpha_{m}} < \alpha_{m+1} \ldots \alpha_{2m}.
\end{equation}
\end{itemize}
\end{lemma}

\begin{proof}[Proof of Theorem~\ref{t11}] Since the set $\uuu \setminus \uu$ is dense in $\uuu$, it is sufficient to show that $\overline{\mathcal{A}} \supset \uuu \setminus \uu$. 
In order to do so, fix $q \in \uuu \setminus \uu$.
Then, according to Theorem ~\ref{t22}, the quasi-greedy $q$-expansion $(\alpha_i)$ satisfies \eqref{23} and 
\eqref{24}. Let $k$ be a positive integer for which equality holds in \eqref{23}, i.e.,  
\begin{equation*}
(\alpha_i)= (\alpha_1 \ldots \alpha_k)^{\infty}.
\end{equation*}
According to Lemma~\ref{l23} there exists an integer $m \geq k$ 
such that \eqref{25} is satisfied. 
Let $N$ be a positive integer  
such that $kN \geq m $ and consider the sequence 
$$(\gamma_i)= (\gamma_i^N)=( \alpha_1 \ldots \alpha_k)^N (\alpha_1 \ldots \alpha_m 
\overline{\alpha_1 \ldots \alpha_m})^{\infty}.$$
For ease of exposition we suppress the dependence of $(\gamma_i)$ on 
$N$. Note that $\gamma_i=\alpha_i$ for $ 1 \leq i \leq m+kN$. In particular, we have 
\begin{equation}\label{27}
\gamma_i=\alpha_i \quad \mbox{for} \quad 1 \leq i \leq 2m.
\end{equation}
Since $(\gamma_i)$ has a periodic tail, the number $q_N$ determined by 
$$1 = \sum_{i=1}^{\infty} \frac{\gamma_i}{q_{N}^i}$$ 
is an algebraic number and $q_N \to q$ as $N \to \infty$.

According to Theorem~\ref{t21} it remains to verify the inequalities \eqref{21} and \eqref{22}. 
First we verify \eqref{21} and \eqref{22} for $j \geq kN$. For those values of $j$ the inequality \eqref{21} for $j+m$ is equivalent to \eqref{22} for $j$ and \eqref{22} for $j+m$ is equivalent to \eqref{21} for $j$. Therefore it suffices to verify the inequalities \eqref{21} and \eqref{22} for $kN \le j < kN +m$. Fix $kN \le j < kN +m$. From \eqref{23}, \eqref{26} and \eqref{27} we have 
\begin{eqnarray*}
\gamma_{j+1} \ldots \gamma_{kN+2m}&=& \alpha_{j-kN +1} \ldots \alpha_m \overline{\alpha_1 \ldots \alpha_m} \\
& < & \alpha_{j-kN +1} \ldots \alpha_m \alpha_{m+1} \ldots \alpha_{2m} \\
&\leq& \alpha_1 \ldots \alpha_{kN + 2m - j}\\
&=& \gamma_1 \ldots \gamma_{kN + 2m - j}
\end{eqnarray*}
and from inequality \eqref{25} we have  
\begin{eqnarray*}
\overline{\gamma_{j+1} \ldots \gamma_{kN+m}} &=& \overline{\alpha_{j-kN+1} \ldots \alpha_m} \\
&<& \alpha_1 \ldots \alpha_{kN+m-j} \\ 
&=& \gamma_1 \ldots \gamma_{kN+m-j}.
\end{eqnarray*}
Now we verify \eqref{21} for $j < kN$. If $m \leq j < kN$, then by \eqref{23} and \eqref{26}, 
\begin{eqnarray*}
\gamma_{j+1} \ldots \gamma_{kN+2m} &<& \alpha_{j+1} \ldots \alpha_{kN+2m} \\
& \leq& \alpha_1 \ldots \alpha_{kN+2m-j} \\
& = & \gamma_1 \ldots \gamma_{kN+2m-j}.
\end{eqnarray*}
If $1 \leq j < m$, then by \eqref{23} and \eqref{25}, 
\begin{eqnarray*} 
\gamma_{j+1} \ldots \gamma_{kN+m+j} &=& \alpha_{j+1} \ldots \alpha_{kN+m} 
\overline{\alpha_1 \ldots \alpha_j}\\
& \leq & \alpha_1 \ldots \alpha_{kN+m-j}\overline{\alpha_1 \ldots \alpha_j}\\
& < & \alpha_1 \ldots \alpha_{kN+m-j} \alpha_{m-j+1} \ldots \alpha_m \\
&=& \gamma_1 \ldots \gamma_{kN+m}.\\
\end{eqnarray*}
Finally, we verify \eqref{22} for $j < kN$. 
Write $j=k \ell +i \, , 0 \leq \ell < N$ and $0 \leq i < k$. If $i=0$, then \eqref{22} follows from the relation 
\begin{equation*}
\overline{\gamma_{j+1}}=\overline{\alpha_1}= 0 < \alpha_1 = \gamma_1.
\end{equation*}
If $ 1 \le i < k$, then applying Lemma~\ref{l23}(ii) we get 
\begin{equation*}
\overline{\alpha_{i+1} \ldots \alpha_{2i}} < \alpha_1 \ldots \alpha_i.
\end{equation*}
Hence
\begin{eqnarray*} 
\overline{\gamma_{j+1} \ldots \gamma_{j+k}} &=& \overline{\alpha_{j+1} \ldots \alpha_{j+k}} \\
&=& \overline{\alpha_{i+1} \ldots \alpha_{i+k}} \\
&<&  \alpha_1 \ldots \alpha_k \\
&=&  \gamma_1 \ldots \gamma_k.\\
\end{eqnarray*}
(In order for the first equality to hold in case $\ell = N-1$, we need the condition $m \geq k$.)
\end{proof}

\begin{remarks} \mbox{}
\begin{itemize}
\item Since the set $\uuu$ is a perfect set and $\uuu \setminus \uu$ is countable, each neighborhood of $q \in \uu$ contains uncountably many elements of $\uu$. Hence the set of transcendental univoque numbers is dense in $\uuu$ as well. 
\item Recently, Allouche, Frougny and Hare (\cite{AHF}) proved that there also exist univoque Pisot numbers. In particular they determined the smallest three univoque Pisot numbers. 
\end{itemize}
\end{remarks}


\begin{thebibliography}{[K$^3$1]}

\bibitem{AC} J.-P. Allouche, M. Cosnard, 
{\em The Komornik--Loreti constant is transcendental}, 
Amer.  Math. Monthly {\bf 107} (2000), no. 5, 448--449.
 
\bibitem{AHF} J.-P. Allouche, C. Frougny, K.G. Hare, {\em On univoque Pisot numbers}, Mathematics of Computation, to appear.

\bibitem{BK} C. Baiocchi, V. Komornik, {\em Quasi-greedy expansions and lazy expansions in non-integer bases}, manuscript. 

\bibitem{EHJ} P. Erd\H os, M. Horv\'ath, I. Jo\'o, 
{\em On the uniqueness of the expansions} $1=\sum q^{-n_i}$, 
Acta Math. Hungar. {\bf 58} (1991), no. 3--4, 333--342. 

\bibitem{EJK1} P. Erd\H os, I. Jo\'o, V. Komornik,   
{\em Characterization of the unique expansions $1=\sum
q^{-n_i}$ and related problems},   Bull. Soc. Math. France  {\bf 118} 
(1990), no. 3, 377--390.

\bibitem{EJK2} P. Erd\H{o}s, I. Jo\'o, V. Komornik,
{\em On the number of $q$-expansions}, Ann.\ Univ.\ Sci.\ Budapest.\
E\"otv\"os Sect.\ Math.\ {\bf 37} (1994), 109--118.


\bibitem{KL1} V. Komornik and P. Loreti,  
{\em Unique developments in non-integer bases},  
Amer.  Math. Monthly {\bf 105} (1998), no. 7, 636--639.

\bibitem{KL2} V. Komornik, P. Loreti, A. Peth\H o,  
{\em The smallest univoque number is not isolated}, 
Publ. Math. Debrecen {\bf 62} (2003), no. 3--4, 429--435.

\bibitem{KL3} V. Komornik, P. Loreti,  
{\em On the  topological structure of univoque sets}, 
J. Number Theory {\bf 122} (2007), 157--183. 

\bibitem{P} W. Parry, 
{\em On the $\beta$-expansion of real numbers}, 
Acta  Math. Hungar. {\bf 11} (1960), 401--416. 

\bibitem{R} A. R\'enyi,
{\em Representations for real numbers and their ergodic properties},
Acta  Math. Hungar. {\bf 8} (1957), 477--493.

\end{thebibliography}
\end{document}